\documentstyle[11pt]{article}
\begin{document}
\def\baselinestretch{1.2}
\newcommand{\cA}{{\cal A}}
\newcommand{\cB}{{\cal B}}
\newcommand{\cC}{{\cal C}}
\newcommand{\cD}{{\cal D}}
\newcommand{\cE}{{\cal E}}
\newcommand{\cF}{{\cal F}}
\newcommand{\cG}{{\cal G}}
\newcommand{\cH}{{\cal H}}
\newcommand{\cI}{{\cal I}}
\newcommand{\cJ}{{\cal J}}
\newcommand{\cK}{{\cal K}}
\newcommand{\cL}{{\cal L}}
\newcommand{\cM}{{\cal M}}
\newcommand{\cN}{{\cal N}}
\newcommand{\cO}{{\cal O}}
\newcommand{\cP}{{\cal P}}
\newcommand{\cQ}{{\cal Q}}
\newcommand{\cR}{{\cal R}}
\newcommand{\cS}{{\cal S}}
\newcommand{\cT}{{\cal T}}
\newcommand{\cU}{{\cal U}}
\newcommand{\cV}{{\cal V}}
\newcommand{\cW}{{\cal W}}
\newcommand{\cX}{{\cal X}}
\newcommand{\cY}{{\cal Y}}
\newcommand{\cZ}{{\cal Z}}

\newcommand{\ST}{{-\!\!\!-\!\!\!-\!\!\!-\!\!\!-\!\!\!
-\!\!\!\!\!\!-\!\!\!\!\!\!-\!\!\!\!\!\!-\!\!\!\!\!\!-\!\!\!\!\!\!\longrightarrow}}
\newcommand{\DT}{{\begin{array}{ll}
\!\!\!\!\mid & {}  \\[-2mm]
\!\!\!\!\mid & {} \\ [-2mm]
\!\!\!\!\mid & {}  \\[-2mm]
\!\!\!\!\mid & {}   \\[-2mm]
\!\!\!\!\!\vee & {}
\end{array}
}}

\begin{center}
{\Large{\bf Boundary of Picard group of a singular curve}}
{\bf Jyotsna Gokhale}\\
Mehta Research Institute\\
10 Kasturba Gandhi Street\\
Allahabad 211002\\
India\\
E-mail : jgokhale@mri.ernet.in
\end{center}
\underline{\bf Introduction:}

Picard group
$ Pic \ C $
of a curve $C$ is a projective variety iff $C$  is non-singular.
The most natural compactification of
$ Pic^0 C \ {\rm is} \ (Pic^0 C)^{=}, \ $
the family of torsion - free sheaves of degree $0$ and rank one on $C$ [MM];
it is irreducible if and only if the embedding
dimension of the $C$ is at most two at every point on $C$ [D], [AIK],
[R], [K]. It is enough to look at one singularity at a time [SK].
In this work
\footnote{ This work was done during '84 - '86  at Brandeis as part of my
doctoral dissertation, and forms the first chapter thereof. Some changes
were made during May '92 before sending this paper for publication, most
notable being the formulation of a functor and proof of its
representability by a projective scheme, with consequent change of
language of the paper. The concept behind this functor, and the use that
has been made of it, is same as in my dissertation.

It was brought to my notice (August '92) that G.Pfister and
J.H.M.Steenbrink have recently ('91) published a paper similar in content
to the second part of my dissertation. I have taken care therefore to
refer to [PS] wherever there is any part common, omitting repetition of
already published proofs.

It is a duty and a pleasure to thank my advisor David Eisenbud for all the
insight and the constructive comments he provided me with, and for the
atmosphere highly conducive to work that he generates around him.}
we study modules over the local ring of an integral curve $C$ with one
singularity $Q$ of embedding dimension more than two
to determine
$ (\overline{Pic^0 C}) $
in
$ (Pic^0 C)^{=}. \hfill{\Diamond}$
\begin{center}
\underline{\bf Notation:}
\end{center}
\begin{tabular}{lll}
$k$                     & {} & an algebraically closed field
                                      of characteristic zero \\
$C$                     & {} & an integral curve over $k$ with
                                 one unibranch singularity $Q$ \\
$C^{\sim}$              & {} & normalisation of $C$ \\
$\pi$                   & {} & the normalisation map \\
$\cO$                   & {} & completion of local ring of $C$ at $Q$ \\
$\cO^{\sim}$            & {} & normalisation of $\cO$ \\
$\cM$                   & {} & maximal ideal of $\cO$ \\
$\cC$                   & {} & the conductor $Hom_{\cO}(\cO^{\sim},\cO)$ of $C$
\\
$\cM_Q$                 & {} & maximal ideal sheaf of $Q$ \\
$\cC_Q$                 & {} & conductor sheaf \\
$\cW^0_C$               & {} & the canonical dualising sheaf of $C$ \\
$\delta$                & {} & $ rk_k (\cO^{\sim} / \cO);$
                                we omit the subscript k normally \\
$\cF$                   & {} & torsion free sheaf of rank one \\
$\cL $                  & {} & a locally free sheaf of rank $1$ \\
$C^D$                   & {} & a partial normalization of $C$ defined
                       (section 1) uniquely for any $ D \subseteq \pi^*(Q) $ \\
$C'$                    & {} & first partial normalization of unibranch $C$
such that \\
    {}                  & {} & $ \ rk_k(\cO'/ \cO) \ = \ 1 \ $
                                 and $ \cC' \neq \cC \ $;
                          i.e.,  $ \ \cO' \ = \ \cO + t^{-1}\cC \ $. \\
$C^i$                   & {} & successive partial normalizations
                                     defined along the same line \\
$\cM^i$                 & {} & maximal ideal of $\cO^i$ \\
$\cC^i$                 & {} & the conductor of $C^i$ \\
$Pic^0 C$               & {} & Picard group of degree $0$ of $C$ \\
$\pi^*$                 & {} & the map $Pic^0 C \longrightarrow Pic^0 C^{\sim}$
                                    of Picard groups corresponding to $\pi$ \\
$K$                     & {} & $ \ ker(\pi^*)$  \\
$(Pic^0 C)^{=}$         & {} & moduli of torsion free sheaves
                                             of degree $0$ and rank one. \\
$\overline{S}$          & {} & Zariski closure of a variety or scheme $S$\\
$(\overline{Pic^0 C})$  & {} & Zariski closure of $Pic^0C$ in $(Pic^0 C)^{=}$\\
$B(\overline{Pic^0 C})$ & {} & boundary of $\overline{Pic^0 C}$, def(0.1) \\
$\chi(\cF)$             & {} &  Euler characteristic of $\cF$ \\
$deg(\cF)$              & {} & $\chi(\cF)-\chi(\cO_C)$ \\
$\Gamma$                & {} & $\Gamma (C) = (k_1, \cdots, k_r)\ {\bf N}$,
                                     the semigroup of orders of generators \\
      {}                & {} & of $\cO$ over $k$, assuming
                                $ 2 < k_1 < k_2 < \cdots < k_r $ \\
$v_0$                   & {} & order of the conductor, i.e.,
                                   $ \cC = t^{v_0}\cO^{\sim} $ \\
$t$                     & {} & any local parameter of
                                   $\cO^{\sim}. \hfill{\Diamond}$\\
\end{tabular}
\footnote{
$ \Gamma^* $
of [PS] is analogous to our $C'$ in case of a unibranch singularity.}
\newpage
\noindent{\underline{\bf  0: Preliminaries}}

We define
$ \ (Pic_{C}^{d})^=(S) \ = \ \{ $
isomorphism classes of torsion free
$ \ \cO_C(S) $
-modules of degree $d$ and rank 1 \}.
$ \ (Pic^{d}_{C})^{=} \ $ is representable and is
represented by a projective scheme
$ \ (Pic^d C)^{=} \ $ \
[D], [AK], which is irreducible [R], [K], iff
$ \ rk(\cM / \cM^2) \leq 2. \ $
We define the boundary of
$ \ \overline{Pic \ C}, \ $ a la [R]:
\medskip\\
{\bf Definition 0.1:} $ \ \cF \in B(\overline{Pic^0 \ C}) \ $
if
$ \ \cF \not \in Pic^0 C, \  {\rm and \ there \ exists} \ \cG, \ {\rm an} \\
\cO_{C \chi \ Spec \ k\left[[t]\right]} \ {\rm -module \ which \ is \
flat \ over} \ k\left[[t]\right], \ $
such that
$ \ (\cG / t \cG) \approx \cF \ $
and the stalk of $\cG$ at the generic point of
$ \ Spec(k\left[[t]\right]) \ $
is a locally free  $\cO_C$-module.
$\hfill{\parallel}$
\medskip\\
In other words, we have a one-parameter family,
$ \ \{ \cL_{\beta} \} \ \subseteq \ Pic^0 C, \ $
with $\cF$ in the closure. Let
$ K = ker (\pi^*) $
where
$ \ \pi^* \ :  \ Pic^0 C \ \longrightarrow \ Pic^0 C^{\sim} \ $
is the map of respective Picard groups corresponding to the normalization map
$ \ \pi : C^{\sim} \longrightarrow C. \ $
Then Rego proves
\medskip \\
{\bf [R,1.2]:}  $ \ \cF \in \ \overline{Pic^0 C} \ \ {\rm iff} \ \
\exists \ \cL \in Pic^0 C \ {\rm such \ that} \ \
\cF \otimes \cL \in \overline{K}. \hfill{\parallel}$
\medskip \\
For
$ d \leq rk(\cO^{\sim}/\cC) $
Rego defines a functor
$$
\begin{array}{lll}
E(\cC, d)(S) & = & \{ {\rm isomorphism \ classes \ of}
\ \cO_S {\rm -modules} \ F_S \mid  \\ {} & {} &
\cC \otimes_k \cO_S \subseteq F_S \subseteq \cO^{\sim} \otimes \cO_S,
{\rm and} \ rk(\cO^{\sim} \otimes \cO_S / F_S)  =  d \}
\end{array}
$$
and shows that it is representable by a projective scheme,
also called
$ E(\cC, d); $
we identify
$ E(\cC, d) $
with its image in
$ (Pic^0 C)^{=}. $
In particular, since
$ K \subseteq E(\cC, \delta), \ $
every boundary point ``defines an element of
$ \ E(\cC, \delta)''\ $ [R, Th.2.3(b)].
So it is enough to consider limits of one - parameter families of free
$\cO$-modules
$ \ \{ F_{\beta} \} \ \subseteq \ E(\cC, \delta) \ $
parameterizsed by
$ \ k\left[[\beta]\right]. \ $
A free module
$ \ \cG \ \in \ E(\cC, \delta)(k(\beta)) \ $
is necessarily of the kind
$ \ \cG \ = \ \partial_{\beta} \cO, \
\partial_{\beta} \ \in \ [\cO \otimes k(\beta)]^* \ $.
We shall write
$ \ \partial_{\beta}\cO \ \longrightarrow \ F \ $
if $F$ is in the closure of such a family.
\medskip\\
\underline{\bf Example 1:}
$ \ \cO = \ k[s^5, s^2t^3, st^4, t^5],
\Rightarrow \ \cO' = \ k[s^3, st^2, st^3], \ $
and
$ \ \cM \ = \ (t^3, t^4, t^5). \ $
Let
$ \ W \ = \ \cO + t\cO \ \approx \ \cW^0_{C,Q} \ $
be the stalk of canonical dualizing sheaf
$ \ \cW^0_C. \ $
We can construct deformations, using
$ \ E(\cC, \delta), \ $
as follows: since
$ \ \cO' \approx t \cO' \ \in \ E(\cC, \delta) \ $
with
$ \ ord_t (t \cO') = 1, \ $
look for
$ \ \partial_{\beta} \in (\cO^{\sim})^{*} \ $
of the form
$ \ \partial_{\beta} = t + \beta \partial'_{\beta}, \
\partial'_{\beta} \in (\cO^{\sim})^{*}. \ $
In fact,
$t + \beta$ will do :
$ \ (t+\beta)\cO = (t + \beta, t^3, t^4, t^5)
\longrightarrow (t, t^3, t^4, t^5) = t \cO'. \ $
Similarly
$ \ \cM \approx t^{-1} \cM \in E(\cC, \delta); \ $
$ \ (t^2 + \beta a_1 t + b^2a_0) \cO  \longrightarrow  t^{-1}\cM \ $
for any
$a_1 \in K$
and any
$ \ a_0 \in k^*. \ $
Up to isomorphism
$ \ \cO, \cO', \cM \ {\rm and} \ W $
the only $\cO$-modules; since for any
$ F \approx W, \ $
necessarily
$ F \not \in E(\cC, \delta), \ $
$W$ has
no possible deformation by free $\cO$-modules.
Therefore
$ \ (Pic^0 C)^{=} \ $
has one component other than
$ \ \overline{Pic^0 C}, \ $
namely, the Zariski closure of
$ \ Pic^2 C \ $-orbit of
$ \ \cW^0_C.  \hfill{\diamond}$
\medskip \\
We will close the section with some remarks assumed and used in
in this work. Since
$ v_0 - \delta $
is used frequently, we will use notation
$ v_0 - \delta = \gamma \ $
for ease of reading.
\smallskip\\
\underline{\bf Remarks:}
\begin{enumerate}
\item[(0.2)] Since
$ \# \{ j \in {\bf N} \mid j \not \in \Gamma \} = \delta, $
so
$ \# \{ j \in \Gamma \mid j < v_0 \} = \gamma. $
\item[(0.3)]
$ j \in \Gamma \Longrightarrow \ v_0 -1-j  \not \in \Gamma \
\forall j \in {\bf Z}. \ $
The converse is true iff $C$ is Gorenstein [HK].
\item[(0.4)] (i)
$ \delta - 1 \leq v_0 \leq 2 \delta $. \\
(ii)
$ \# \{ j \in {\bf N} \mid j, v_0 -1-j \not \in \Gamma \} = 2\delta -v_0$.
\item[(0.5)] If $F$ is an $\cO$-module such that
$ \cC \subseteq F \subseteq \cO^{\sim} \ $,
and
$ End(F) $ does not contain
$ t^{-1} \cO \ $,
then
$ \cO \subset F \subseteq \cO^{\sim}, \ $
and
$ v_0 - \delta \leq rk(F / \cC) < \delta. $
Further, if
$ rk(F / \cC) = \delta $
then
$ F \approx \cW^{0}_{C,Q} \ $ necessarily. Conversely, for
$ j : v_0 - \delta \leq \ j \leq \delta, $
there exists an $F$ such that
$ \cC \subseteq F \subseteq \cO^{\sim}, \ rk(F / \cC) = j, \ $
and
$ End(F) $
does not contain
$ t^{-1} \cC. $
\item[(0.6)] The following are equivalent: \\
(i) $ v_0 = k_{\gamma-1} + k_1. $ \\
(ii) $ k_i = ik_1 \ \forall \ i \leq \gamma, \ k_i \in \Gamma. $ \\
(iii) $ \ rk(\cM / \cM^2 + t\cC) = 1.$
\item[(0.7)]
$ E(\cC, \delta) \ {\rm is \ irreducible}
\ \Longleftrightarrow \ \overline{K} \ = \ E(\cC, \delta). \ $
\item[(0.8)]
$ \ \cM_Q \ \in \ \overline{Pic \ C} : \ $
the family
$ \ \{ \cO_C (-P) \mid \ P \in C, \ P \ \neq \ Q \} \ $
of locally free sheaves of rank one and degree -1 has $\cM_Q$  in the closure.
\item[(0.9)] If
$ \ \cF \ \in \ B(\overline{Pic^0 C}) \ $
then
$ Hom_{\cO_c} (\cF, \cO_c) \ = \ \cF^v \ \in \ B(\overline{Pic^0 C}). $
This follows from the definition of
$ \ B(\overline{Pic^0 C}). \hfill{\Diamond}$
\end{enumerate}
\noindent{\underline{\bf  1: General Singularity}}

Let $C$ be an integral curve singular at one point
$ Q \in C. $
Then
$ \cO_{C^{\sim}, Q} = \cO^{\sim}, $
the normalization of
$ \cO_{C, Q} = \cO $
in the function field of $C$, is a semi-local domain, with maximal ideals
$ \cM_1, \cM_2, \cdots \cM_k $
corresponding to the points of
$ Supp(\pi^*Q) = \{Q_i \mid i = 1, \cdots , k \} \subseteq C^{\sim}. \ $
Let
$ \pi^*Q = \sum \{ \mu_i Q_i \mid i = 1, \cdots, k \}, $
and let
$t_P$
denote a local parameter of
$ \cO^{\sim} $ at $P$.
Let  $\cC$ be the conductor
$ Hom_{\cO}(\cO^{\sim}, \cO)$.  Since $\cC$ is an ideal of
$ \cO^{\sim} $
as well as $Q$,
$ \ \cC = \cM_1^{v_1} \cdots \cM_k^{v_k} \ $
for some
$ \{ v_i \geq \mu_i \mid i = 1, \cdots , k \}. $
Let $\cM$ be the maximal ideal of $\cO$.  We have
$ \ \cM \subseteq \cM_1^{\mu_1} \cdots \cM_k^{\mu_k}, \ $
and equality holds if and only if
$ \ \cC = \cM. $
Moreover,
$ \ \exists \ g \in \cM $
such that for any
$ P \in \pi^* Q, \ (g)_P = u_P (t_{P}^{\mu_P}) $
for some
$ u_P \in \ {\cO}_{P}^{\sim *}. \ $
Also, if
$ \mu = \sum \mu_i $
and
$ v = \sum v_i, $
then
$ \mu-1 \leq \delta = rk_k(\cO^{\sim} / \cO) \leq v-1, $
and
$ \delta = \mu-1 $
or
$ \delta = v-1 $
iff
$ \cM = \cC. $

We define
$ \ \cO^{Q_i} = \cO + \cM_{1}^{v_1} \cdots
\cM_{i-1}^{v_{i-1}} \cM_i^{v_{i}-1}\cM_{i+1}^{v_{i+1}} \cdots \cM_{k}^{v_k} \ $
for any
$ Q_i \in \pi^*Q, $
and
$ \cO^D $
for
$ D \subseteq \pi^*Q $
successively (by normalizing one point at a time).
There exist partial normalizations $C^P$ of $C$, for any
$ P \in \pi^*Q, $
with
$ \cO_{C^P, Q} = \cO^P, $
$ \cO_{C^P} \mid_{C-Q} = \cO_C \mid_{C-Q}. $
Therefore
$ rk_k (\cO_{C^P} / \cO_C) = 1 $
and
$ \cC^P \neq \cC. $
Let
$ \pi^P : C^P \longrightarrow C $
be the partial normalization map corresponding to the natural inclusion
$ \cO \hookrightarrow \cO^P. $
\medskip\\
{\bf Lemma 1.0:} The following are equivalent \\
(i) $ \cM \approx \cO^D $ for some $ D \subseteq \pi^*Q $ \\
(ii) $ rk(\cM / \cM^2 + \cC) = 1 $ \\
(iii) $ \cM = g \cO + \cC $ for some $ g \in \cM. $
\smallskip\\
{\bf Proof:}  Given an isomorphism
$ f : \cO^D \longrightarrow \cM, $ let
$g = f(1)$,
so that
$ m = g \cO + g \cC^D. $
This gives
$ g \cC^D \subseteq \cM, \Longrightarrow \ g \cC^D \subseteq \ \cC, $
by definition of $\cC$. Converse is obvious: if
$ \cM = g \cO + \cC $
then
$ g^{-1} \cC \subseteq \cO^{\sim}, $
and
$ \cM  = g \cO^D $
for
$ \cO^D = \cO + g^{-1} \cC. \hfill{\parallel}$
\medskip\\
{\bf Theorem 1.1:}  $ \ B(\overline {Pic^0 C}) \ \  \subseteq
\ \ \bigcup \ \ \{\pi^{P}_{*} (Pic^{-1} C^{P^{=}}) \ \mid P \in \pi^*(Q) \}. $
\smallskip\\
{\bf Proof:}  We need to prove that if
$ \ \cF \in \overline{Pic^0 C} \ $
such that
$ \cF_Q $
is not isomorphic to
$ \cO_Q $
then there exists a point
$ P \in Supp \{ \pi^* Q \} $
such that
$ \cF_Q $
is an
$ \cO^P $-module. If
$ \cF \in \overline {Pic^0 C}, \ $
then there exists
$ F \approx \cF_Q $
such that
$ \ \cC_Q \subseteq F \subseteq \cO^{\sim} \ $
and
$ rk(F / \cC) = \gamma. $
So
$ \cF_Q $
is not isomorphic to
$ \cO \Longleftrightarrow F $
contains no unit of
$ \cO^{\sim} \Longleftrightarrow  F \subseteq \cM_P $
for some
$ P \in \pi^* Q $
where
$ \cM_P $
the maximal ideal of  $P$  in
$ \cO^{\sim} \Longleftrightarrow F $
is an
$ \cO^P $-module for some
$ P \in \pi^*Q. \hfill{\parallel} $
\medskip\\
{\bf Theorem 1.2:}
$ \ \bigcup \
\{ \pi^{P}_{*} \ (\overline {Pic^{-1} C^P}) \mid \ P \in \pi^*(Q)\} \
\subseteq \ B(\overline{Pic^0 C}). $
\smallskip\\
{\bf Proof:}  Enough to prove that
$ \pi^{P}_{*} (Pic^{-1} C^P) \subseteq \overline{Pic^0 C} \ \forall i, $
or, equivalently, given
$ P \in \pi^* Q, $
there exists
$ \cF \in \overline{Pic^0 C} $
such that
$ \cF_Q \approx  \cO^P. $
Let  $t_P$  be a local parameter of
$ \cO^{\sim} $
at $P$ and let
$ \partial_{\beta} = t_P - \beta. $
This will do:
$ (t_P - \beta) \cO \longrightarrow t_P \cO^P \approx \cO^P. \hfill{\parallel}$
\medskip\\
\underline{\bf Example 2:} \ If
$ \ \cO = k \left[[t^3, t^7, t^8]\right] \ $
then
$ \ \cO' = k \left[[t^3, t^5, t^7]\right] \ $
has a deformation by $\cO$-modules, but
$ A = k \left[[t^3, t^4]\right] \ $
does not, (nor is
$\cO'$
an $A$-module,) even though $A$ is a partial normalization of singularity
one less than that of
$\cO$
just as
$\cO'$
is.
So the definition of
$ C^P $
is relevant to determination of
$ B(\overline{Pic \ C}). $
$\hfill{\diamond}$
\medskip\\
{\bf Theorem 1.3:} $ B(\overline{Pic^0 C}) = \bigcup \ \{
\pi^D_{*} (Pic^{-i} C^D) \mid \ D \subseteq \pi^* Q, \ {\rm deg} \ D = i \} $
if and only if
$ \ \cM = \cC. $
\smallskip\\
{\bf Proof:}  If $ \cM = \cC $,
then
$ \delta = v_0 -1, $
and therefore
$ E(\cC, \delta) = E(\cC, v_0 -1). \ $
Since
$ \overline{Pic^0 C} \subseteq E(\cC, \delta), \  \cF \in
\overline{Pic^0 C} \Longrightarrow \exists F \approx \cF_Q $
such that
$ \cC \subseteq F \subseteq \cO^{\sim} $
and
$ rk(F / \cC) = 1. $
Let
$ a \in F, \ a \not \in \cC $:
then
$ a \cO^{\sim} = \cM_1^{i_1} \cdots \cM_k^{i_k}, $
and
$ F = a \cO + \cC \approx \cO^D$ for $D = \sum i_j Q_j. $
This implies
$ \cF \in \pi^{D}_{*}(Pic^{-i} C^D). $

Conversely, suppose
$ \cM \neq \cC. $
If $\cM$ is not isomorphic
to $\cO^D$ for any
$ D \subseteq \pi^*Q, $
then the sheaf
$ \cM_Q $
is a counterexample, since
$ \cM_Q $
has a deformation by locally free sheaves.  If on the other hand
$ \cM \approx \cO^D $
for some
$ D \subseteq \pi^* Q $
then
$ \cM = g \cO + \cC $
for some
$ g \in \cO. $
We will show the existence of $f$ such that
(i) $ \cM + f \cO $
has a deformation by free $\cO$-modules, and
(ii) $ \cM + f\cO $
is not isomorphic to $\cO^D$ for any
$ D \subseteq \pi^* Q. $

Let
$ r = rk (\cM / \cC). $
For
$ P \in \pi^*Q $
let $t_P$  be a local parameter at $P$.  If
$ g^r t_P \cO + \cM $
is not isomorphic to $\cO^P$, we set
$ f = t_P g^r. $
If
$ g^{r} t_{P} \cO + \cM \approx \cO^P $
then
$ t_P^{-1} g^r \cO + \cM $
is not isomorphic to any
$ \cO^D, $
so
$ t_P^{-1} g^r= f $
will do.  In either case
$ F = f \cO + \cM $
is not isomorphic to any
$ \cO^D, $
and
$ (f + \beta) \cO = (f + \beta , g) \cO \supseteq \cM, $
so
$ (f+\beta) \cO \longrightarrow F = f \cO + \cM $
and hence the theorem.
$\hfill{\parallel}$
\medskip \\
The first examples of a reducible
$ \ (Pic^0 C)^{=} \ $
are provided by
$ \delta = 2. $
\medskip\\
{\bf Lemma 1.4:}  If
$ \ \delta = 2, \ {\rm then}
\ rk(\cM / \cM^2) = 2 \ \Longleftrightarrow \ \cM \neq \cC. $
\smallskip\\
{\bf Proof:}  If
$ \cM = \cC, $
then
$ \ rk(\cM / \cM^2) = rk(\cC / \cC^2) = rk(\cO^{\sim} / \cC)
= rk (\cO^{\sim}/\cM) = 3. $
If not, we have
$ rk(\cO^{\sim} / \cM) = \delta + 1 > \mu, $
so that
$ \mu < 3. $
But
$ rk(\cO^{\sim} / \cC) > rk (\cO^{\sim} / \cM) = 3,
\ \Longrightarrow \ \sum_1^k v_i > 3. $
Now,
$ \cO^{\sim} $
has at most two maximal ideals
$ \cM_1 $
and
$ \cM_2, $
and
$ \cM \subset \cM_1 \cM_2 $
as an $\cO$-submodule; if
$ \cM \subseteq \cM_1^2 \cM_2, $
then since
$ rk(\cO^{\sim} / \cM_1^2 \cM_2) = 3 $
we must have
$ \cM = \cM_1^2 \cM_2, $
and therefore
$ \cC = \cM, $
a contradiction.  Thus
$ \cM $
is not contained in
$ \cM_i^2 \cM_j. $
For the same reason
$ \cC \neq \cM_i^2 \cM_j, \ \Longleftrightarrow \
\cC \subseteq \cM_1^2 \cM_2 \bigcap \cM_1 \cM_2^2, \ $
since
$ v_1, v_2 > 1. $
Therefore
$ \ \exists a \in \cM $
such that
$ a = a_1 a_2, \ a_i \in \cM_i - \cM_i^2, \ a \not \in \cM_i^2 \cM_j; $
i.e.,
$ a_i $
must generate
$ \cM_i $
locally.  But then
$ rk(\cM / \cM^2) = 2, $
since
$ rk(\cM / \cM^2) \leq rk(\cC / \cM \cC) \leq rk(\cC / a \cC) = 2.
\hfill{\parallel}$
\bigskip\\
\underline{\bf 2: $Filt(\cC, \delta)$:}

For rest of this work $C$ will be assumed to have a unibranch singularity
$ Q; $ i.e.,
$\pi^*Q $
is supported at one point
$ Q \in C^{\sim}. $
(See Notation).
$ \cO $
has a unique filtration
$ \cO \supseteq I_1 = \cM \supseteq
\cdots \supseteq  I_{v_0 - \delta - 1} \supseteq \cC $
such that
$ rk(I_j / I_{j+1}) = 1 $
and
$ {\rm ord}_t I_j< {\rm ord}_t I_{j+1} $
for all
$ j = 1, \cdots \gamma -1 $ .
Since
$ I_j = (\cO^j)^v, $
from (0.9)
\medskip\\
{\bf Lemma 2.0:}  If
$ \cF \in \ (Pic^0 C)^{=} $
such that
$ \cF_Q  \approx I_j, $
then
$ \cF \in B(\overline{Pic^0 C}).  \hfill{\parallel}$
\medskip\\
An
$ \cO- {\rm module} \ F $
has a unique filtration analogous to that of $\cO$ given by
$ F \supseteq F_1 \supseteq \cdots  $
such that
$ F_{i+1} \subseteq tF_i \cO^{\sim} $
and
$ rk(F_i/F_{i+1}) = 1.$
If moreover
$ F \in  E(C, \delta -r-1) \ $
then
$ F_{r+1} = \cC, $
and if  $F$  has a deformation by free  $ \cO $ -modules then
$ F_i \subseteq t^{k_i} \cO^{\sim} $
for all
$ i = 0 \cdots \gamma -1.$
This motivates the following construction :
\newpage
\noindent{\bf Definition 2.1:}
$ \ Filt^{\delta}_{\cC} (S) = \{ F_S \in E(\cC, \delta)(S) \mid
(F_{S})_j  \subseteq I_j \cO^{\sim} \ \ \forall j \geq 0 \}, $
where
$ (F_{S})_j $
is the
$ j $th  module in the unique filtration
$ F_S \supset (F_{S})_1 \supset \cdots \supset \ \cC \otimes \cO_S $
such that
$ ord_t (F_{S})_j < ord_t (F_{S})_{j+1} \ \forall \ j, \  $
and
$ (F_{S})_j / (F_{S})_{j+1} $
is a rank one locally free
$ \cO_S {\rm -module}.  \hfill{\parallel} $
\smallskip\\
Next we define an extension of
$ E(C, \delta) $
and use it to show that
$ Filt^{\delta}_{\cC} $
is representable by a projective scheme
$ Filt \ (\cC, \delta) = \{F \in E(\cC, \delta) \mid F_j
\subseteq I_j \cO^{\sim}, \  0 \leq j \leq \gamma-1\}. $
Let
$ \widetilde{I} \  $
be any ideal of
$ \ \widetilde{\cO}, \ $
with
$ \cC \subseteq \widetilde{I}, \ $
and let
$ d \leq rk(\widetilde{I} / \cC). $
\medskip\\
{\bf Definition 2.2:}
$ E(\cC, \widetilde{I}, d)(S) = \{\cO \otimes_k \cO_S {\rm -modules} \ F_S
\mid \ \cC \otimes \cO_S \subseteq F_S \\
\subseteq \widetilde{I} \otimes_k \cO_S \
{\rm and} \
\widetilde{I} \otimes \cO_S / F_S \
{\rm is \ a \ locally \ free \ } \cO_S - {\rm module \ of \  rank} \ d \}.
\hfill{\parallel} $
\smallskip \\
Rego's proof of representability of
$ E(\cC, d) $   [R, (2.,1)]
carries over with a minor change -i.e., replacing
$ \widetilde{\cO} $  by
$ \widetilde{I} \ - \ $
so we have
\medskip \\
{\bf Lemma 2.3:}
$ E(\cC, \widetilde{I}, d) $
is representable by a projective scheme.
$ \hfill{\parallel} $
\medskip \\
We shall denote the scheme representing
$ E(\cC, \widetilde{I}, d) $
by the same. For any
$ F \in E(\cC, \delta), $
necessarily
$ F_j \in E(\cC, \overline{I}_j, k_{\gamma-j} - \gamma + j) $
where
$ \overline{I}_j = t^{v_0- k_{\gamma-j}} \widetilde{\cO}. $
Also,
$ t^{k_j} \widetilde{\cO} = \widetilde{I}_j \subseteq \overline{I}_j. $
So for
$ j \leq \gamma -1 $
we have morphisms:
$$
\begin{array}{lllll}
{} & {} & E(\cC, \delta) & {}  & F \\
{} & {} & \varphi_j \downarrow & {} & \downarrow \\
E(\cC, \widetilde{I}_j, d_j) &  \longrightarrow &
E(\cC, \overline{I}_j, \overline{d}_j)& {}  & F_j  \\[-2mm]
{} & id & {} & {} & {}
\end{array}
$$
which are proper since
$ E(\cC, \cI, d) $
are projective. Now
$$ F\in Filt(\cC, \delta) \ \Longleftrightarrow  \
F_j \in E(\cC, t^{k_j} \widetilde{\cO}, \delta+ j - k_j)
\ \ \forall \ \ j \leq \gamma-1, $$
so
$$ Filt(\cC, \delta) \ = \ {\displaystyle{\bigcap_{j=0}^{\gamma}}}
\varphi_j^{-1} (E(\cC, t^{k_j} \widetilde{\cO}, \delta + j-k_j)). $$
This proves
\medskip\\
{\bf Theorem 2.4:} $ \ Filt_{\cC}^{\delta} \ $
is representable by a projective scheme.
$\hfill{\parallel}$
\medskip \\
There is an obvious map
$  Filt^{\delta}_{\cC} \longrightarrow  Pic^{0^{=}}_{C}, $
and the image in
$  (Pic^0 C)^{=}  $
is
$ Filt(\cC, \delta). $
Since
$ \cO \in \ Filt (\cC, \delta), $
we have
$ K \subseteq  \ Filt (\cC, \delta) \subseteq E(\cC, \delta). $
We shall skip the proof of the following:
\medskip\\
{\bf Lemma  2.5.1:}  If
$ rk(\cM / \cM^2 + \cC) = 1 \  $
and
$  \ v_0 \geq k_{\gamma-1} + k_1 - 1, $
then
$ E(\cC, \delta) = Filt (\cC, \delta). $
\medskip\\
{\bf Lemma 2.5.2:} If  $ E(\cC, \delta) = Filt (\cC, \delta), $
then one of the following conditions hold good:\\
(i) If
$ rk(\cM/\cM^2 + \cC) = 1 $  then
$ v_0 \geq k_{\gamma-1} + k_1 -1. $ \\
(ii) If
$ v_0 < k_{\gamma-1} + k_1-1 $
then
$ v_0 - k_1 \leq k_{\gamma-2} $
and
$ k_{\gamma-1} \leq v_0 - k_{\gamma-1} + k_{\gamma-2}; $
moreover for any
$ j \in \{2, \cdots, \gamma-2 \} $
such that
$ v_0 -k_1 < k_j, $
it follows that
$ k_j \leq v_0 - k_{\gamma-1} + k_{j-1}. \hfill{\parallel}$
\medskip\\
In general
$ \overline{K} = Filt (\cC, \delta) $
is not true. Consider
\medskip\\
\underline{\bf Example 3:}
$ \cO = k \left[[t^{k_1}, t^{k_1+2} ]\right] + t^{k_1+4} \cO^{\sim} $
and
$ F = (t^4, t^{k_1+1}, t^{k_1+2}) + \cC.$
If there is an
$ \partial_{\beta} \in \cO^{\sim{*}} $
giving a deformation of  $ F $ , then necessarily
$ \partial_{\beta} = t^4 + \beta \partial'_{\beta} $  and since
$ F_1 = (t^{k_1 +1}, t^{k_1+2})+\cC, \partial'_{\beta}
= t + \beta \partial''_{\beta}. $
This contradicts with
$ {\rm ord}_t F_2 = k_1+2. \hfill{\diamond}$
\medskip\\
{\bf Lemma 2.6:}  If
$ rk(\cM / \cM^2 + \cC) = 1, $
then
$ Filt(\cC, \delta) = \overline{K}. $
\smallskip\\
{\bf Proof:}  Let
$ {\rm ord}_t F = h. $
If
$ rk((t^h \cO + \cC)/\cC) = \gamma, $
then
$ F = t^h \cO + \cC \approx t^h \cO^h $
and  $ F $  has a deformation. If
$ rk((t^h \cO + \cC)/\cC) < \gamma, $
set
$ r = rk(F/t^h \cO+\cC). $
We may assume
$ F = (t^h, f_1, \cdots, f_r) \cO + \cC, $
where  $ f_i $  have increasing orders,
$ v_0-(\gamma -r) k_1 \leq h < v_0-(\gamma - r-1)k_1 $
by definition of  $ r, $  and for
$ \ i = 1, \cdots, r, \
rk[(f_i, \cdots, f_r) + \cC / \cC] \
= \ r - i + 1, \ {\rm so \ that} \
{\rm ord}_t f_i \geq (\gamma - r-1 + i) k_1 . $
Using the condition on
$ {\rm ord}_tf_i, $
we construct
$ \partial_{\beta} $
such that
$ \partial_{\beta} \cO = (\partial_{\beta},
f_1 + \beta f'_1, f_2 + \beta f'_2, \cdots, f_r + \beta f'_r) $
for some
$ \{f'_1. f'_2. \cdots , f'_r \} \subseteq \cO^{\sim}, $
so that
$ \partial_{\beta} \cO \longrightarrow F. $
Let
$ \partial_{\beta}
= t^h + \beta g^{-(\gamma-r)} f_1 + \beta^2 g^{-(\gamma - r - 1)} f_2 +
\cdots + \beta^{g^ -(\gamma-1)} f_r + \beta^{r+1}. $
It is easy to check that we have constructed the requisite
$ \partial_{\beta}.\hfill{\parallel} $
\medskip \\
{\bf Corollary 2.6.1:} If  $ rk(\cM/\cM^2 + \cC) = 1 $  and
$ v_0 \geq k_{\gamma-1} + k_1-1, $
then
$ E(\cC, \delta) =  \overline{K}. $
\medskip\\
\underline{\bf Example 4:} The converse of (2.6) is not true : for
$ \cO =  k\left[[s^6, s^2 t^4, st^5, t^6]\right], $
every  $ \cO $ -module
$ F \in \  Filt\ (\cC, \delta) $
has a deformation by free $\cO$-modules.  This can be seen explicitly : if
$ ord_t F = 1 $
then
$ F = t\cO' $
necessarily; if
$ ord_t F = 2 $
and
$ F \neq t^2 \cO'' $
then
$ F = (t^2, ut^4)+C $
for
$ u \in \cO^{\sim{*}}, $
and
$ \partial_{\beta} = t^2 + \beta u $
will do; if
$ ord_tF = 3 $
then  $ F $   is in the closure of the Picard orbit of
$ (t^3, t^4, t^5)+\cC = t^{-1} \cM, $
and
$ ord_tF = r, \Longleftrightarrow F = t^4 \cO^{\sim}. $
Hence
$ Filt(\cC, \delta) = \overline{K}. \hfill{\diamond}$
\medskip\\
As example(2) shows,
$ Filt(\cC, \delta) = \overline{K} $
is not true in general even in case
$ rk(\cM / \cM^2 + \cC) = 2, \ $
although it holds for
$ rk(\cM / \cM^2 + \cC) = 3 $
in example(3).
However, using (2.5) and (2.6) we have
\medskip\\
{\bf Theorem 2.7:}  If
$ E(\cC, \delta) = \overline{K} $
then one of the conditions of (2.5.2) holds good. Moreover
\begin{enumerate}
\item[(i)]  in case 2.5.2(i), if
$ v_0 =k_{\gamma-1} + k_1-1, rk(\cM/\cM^2 + \cC) \neq 1 $
then
$ k_1 \leq 4, $
\item[] and
\item[(ii)]  in case of 2.5.2(ii), if
$ v_0 - k_1 = k_{\gamma-2} $  then
$ k_1 = 2n + 3 \geq 5 $  and
$ v_0=  k_{\gamma-1}+ k_1-2 = k_{\gamma-3} + 2k_1 -2. $
\end{enumerate}
\newpage
\noindent{\bf Proof:} If
$ E(\cC, \delta) = \overline{K}, $
then trivially
$ E(\cC, \delta) = Filt (\cC, \delta), $
i.e., hypothesis of (2.5.2) holds; so either 2.5.2(i) or 2.5.2(ii) must hold.

To show (i), if
$ v_0 = k_{\gamma-1} + k_1 -1, k_{\gamma-2} + k_1 > k_{\gamma-1} $
and
$ k_1 > 4 $
then
$ F \in E(\cC, \delta) - \overline{K} $
for
$ F = (t^{k_1+4}, t^{v_0 - k_1},
t^{v_0 - k_1+2}, t^{v_0 - k_1 +4}) \cO + \cC. $

To show (ii), Suppose conditions of 2.5.2(ii) hold, so that we have
$ v_0 < k_{\gamma-1} + k_1 -1 $  and
$ k_{\gamma-1} \leq v_0 - k_{\gamma-1} + k_{\gamma-2}; $
and suppose that
$ v_0 - k_1 = k_{\gamma-2}. $

We set
$ h = v_0 - k_{\gamma -1}, F = (t^h, t^{h'})\cO + \cC $
with
$ v_0 - k_1 \leq h' < k_{\gamma-1} $
chosen so that
$ F \in E (\cC, \delta) $
and
$ ord_t F_{\gamma-2} = h'. $
If  $ F $  has a deformation
$ \{a_{\beta} \cO\} $
by free  $ \cO $ -modules then necessarily
$ a_{\beta} = t^h + \beta t^{h+k_{\gamma-2} - k_{\gamma-1}} +
\beta^2, \ \ a_{\beta} \in (\cO^{\sim} \otimes k \left[[ \beta]\right])^* $
in order to have
$ a_{\beta} \longrightarrow  t^h $
and
$ \{a_{\beta} I_{\gamma-1}\} \longrightarrow t^h $
and
$ \{a_{\beta} I_{\gamma-1}\} \longrightarrow F_{\gamma-1}. $
Since we need
$ \{a_{\beta} I_{\gamma-2}\} \longrightarrow F_{\gamma-2} $
either
$ h + 2k_{\gamma-2} - k_{\gamma-1} = h' $
or
$ h+ 2 k_{\gamma-2} - k_{\gamma-1} \geq k_{\gamma-1}, $
so that
$ a_{\beta} $  can be modified further.

In the second case let
$ h_j = h + j(k_{\gamma-2} - k_{\gamma-1},
\ n  =  {\rm max} \{j \mid h_j > 0\} $
and
$ a_{\beta j} = t^h + \beta^{h_1} + \cdots + \beta^j t^{h_j} + \beta^{j+1} $ ,
for
$ j \leq n. $
If
$ \{a_{\beta} \cO\} \longrightarrow F $
then
$ a_{\beta} = a_{\beta j} $
for some  $ j. $  But
$ \{ a_{ \beta j} I_{ \gamma -2 } \} \longrightarrow G_j, $
$ ord_tG_j = \ {\rm min} \ \{k_{\gamma-2} + h_j, k_{\gamma-1} \}. $
Since
$ \{ h_j \} $
are not consecutive we may choose  $ h' $  to avoid them so that
$ F \not \in \overline{K}, $
unless
$ \{ v_0 - k_1, v_0 - k_1 +1, \cdots, k_{ \gamma-1 } \} \subseteq
\{ h+ k_{ \gamma-3 }, h_j, k_{ \gamma-1 } \mid j \leq n \} $
which is only possible for
$ h + k_{\gamma-3} = v_0 - k_{1}, \ k_{\gamma-1} = k_{\gamma-2} +2 $  and
$ h_n = v_0 - k_1 + 1, $
i.e., the conclusion of the statement;
$ k_1 = 2 n +3 \geq 5 $
in this case. Since the case
$ h' = h + 2k_{\gamma-2} - k_{\gamma-1} $  is
included in the discussion above, we are through.
$ \hfill{\parallel} $
\medskip\\
{\underline{\bf Example 5:}} If
$ \cO = k \left[ [ t^4, t^{13}, t^{18}, t^{19} ] \right] $
and
$ \Gamma(F) = \{ 5,9,11, 13, 14 \} \bigcup  \Gamma(\cC) $
then
$ F \not \in \overline{K}. $
So for
$ k_1 \leq 4, $
it is not necessarily true that if
$ \ v_0 = k_{\gamma-1} + k_1 -1 \ $
and
$ \ rk(\cM / \cM^2 + \cC) \neq 1 \ $
then
$ \ E(\cC, \delta) = \overline{K}. \hfill{\diamond}$
\medskip\\
{\bf Lemma 2.8:} If  $ v_0 = k_{\gamma-1} + k_{1}-1 $  and
$ E(\cC, \delta) = \overline{K} \ $
then
$ rk(\cM/\cM^2 + \cC) = 1 $
unless either \
(i)  $ k_1 \mid v_0, k_1 \leq 3 $
\ or \
(ii) \  $ k_1 \mid v_0 + 1, k_1 \leq 4 $
\smallskip \\
{\bf Proof:}  If
$ v_0 = k_{\gamma-1} + k_1 -1 $
and
$ rk(\cM + \cM^2 + \cC) \neq 1 $
then
$ k_{\gamma-2} = v_0 - k_1 $
and
$ k_{i+1} - k_i \in \{ 1, k_1 -1, k_1 \} \forall i. $
Also, not more than two
$k_i$
can be consecutive
$ \forall i < \gamma. $
Therefore, either
$ k_1 \mid k_{\gamma-2} $
or
$ k_1 \mid k_{\gamma-1}, $
and so either \\
(i)  $ \Gamma = \{ k_1, \cdots,
nk_1, nk_1 + 1, \cdots, mk_1, mk_1 + 1; (m+1) k_1 = v_0, \cdots \}, $ \\
\ \ \ \ or \\
(ii)  $ \Gamma = \{ k_1, \cdots, (n-1) k_1, nk_1 -1, nk_1,
\cdots, mk_1 -1, mk_1; (m+1) k_1 -1 = v_0 \cdots \}. $  \\
Set
$ \ F = (t^{k_1}, t^{nk_1+2}, t^{mk_1+3}) \cO+\cC. $
Then
$ F \in Filt (\cC, \delta) $
but
$ F \not \in \overline{K} $
unless
$ k_1 \leq 3 $  in case (i) or
$ k_1 \leq 4 $  in case (ii).
$ \hfill{\parallel} $
\bigskip\\
\underline{\bf 3: $ B(\overline{Pic^0 C} ) $}

We would like to explore the question of when the inequalities in
(1.1) and (1.2) become equalities. The first, namely,
$ B(\overline{Pic^0 C}) = \pi'_{*}(\overline{Pic^{-1} C'}) $
is true for
$ \cM = \cC $,
but also for a much larger category of curves. It ties in with
(2.6.1) and (2.7) to give a better answer to the previous topic. Let
$ \cF \in \overline{Pic^0 C}, \ \cF \not \in Pic^0 C, $
and
$ F \approx \cF_{Q} $
such that
$ F \in E(\cC, \delta). $
If
$ \cF \in \overline{Pic^{-1} C'}, \ $
then
$ \ \exists \ F' \approx F \ $
such that
$ \ F' \in E(\cC', \delta-1). $
Since
$ E(\cC', \delta -1) \subseteq E(\cC, \delta-1), $
necessarily
$ \ tF' = F, \ $
or
$ \ F \supseteq t \cC', \ $
for some local parameter $ t $  at  $ Q $ .
\medskip\\
{\bf Lemma 3.1:}  If
$ \ B(\overline{Pic^0 C}) =
\pi'_{*} (\overline{Pic^{-1} C'}), \ $
then
$ \ rk_k(\cC'/ \cC) = 1 \ $
unless
$ k_1 = 2. $
\smallskip\\
{\bf Proof:}  Suppose
$ k_1 \neq 2. $
If $C$ is Gorenstein, then
$ v_1 = v_0 - k_1. $
Take
$ F = t^{-1} \cM. $
Since
$ t^{-1} \cM \in E(\cC, \delta) $
and the sheaf
$ \cM_Q $
is in the completion of the family
$ \{ \cO_C(-P) \mid P \in C, P \neq Q \} $
we know that an
$ \partial_{\beta} $
giving a deformation of $F$  exists.
On the other hand
$ \cC' $
is not contained in
$ t^{-2} \cM, $
so $F$ has no deformation by free
$ \cO' $
-modules.

If $C$ is not Gorenstein, this argument may not work for
$ v_1 = v_0-2. $
However, we can construct a counterexample by using the fact that
$ 2 \delta-v_0 = \ \#  \{ j \mid j, v_0 -1 - j \not \in \Gamma \} > 0 $
(which, therefore, works only for non-Gorenstein curves,) as follows :
Set
$ \partial_{\beta} = a_0 + a_1 t + \cdots  + a_{v_0-1} t^{v_0-1}, $
where
$ a_i = 0 $  if  $ v_0 -1 - i \in \Gamma, $
and
$ a_i = \beta^{v_0 -1-i} $
otherwise. Let
$ \partial_{\beta} \cO \longrightarrow F $ .
Since $C$ is not Gorenstein,
$ F \neq \cO. $
On ther other hand  $ F $  does not contain
$ t^{-1}C $
by construction, so
$ t^{-1} F \supseteq \cC' $
iff
$ \cC' = t^{-1} \cC.\hfill{\parallel} $
\medskip \\
The converse is not true :
$ rk_k(\cC' / \cC) = 1 $
is a necessary but not sufficient condition for
$ B(\overline{Pic^0 C}) = \pi'_{*} (\overline{Pic^{-1} C'}). $
\medskip\\
\underline{\bf Example 6:} Let
$ \cO $
be given by
$ \cO = k\left[[t^{k_1}]\right] + t^{2k_1-1} \cO^{\sim}, $
and let
$ \partial_{\beta} = t^{k_1 -1} + \beta. $
Then
$ F = (t^{k_1-1}, t^{k_1})\cO + C, \cO'
= k\left[[t^{k_1}]\right] + t^{2k_1-2} \cO^{\sim}, $
and
$ t^{-1} F = (t^{k_1 -2}, t^{k_1-1}) \cO' + \cC'. $
Since
$ t^{-1} F \not\in \ Filt \ (\cC, \delta-1), $
therefore
$ t^{-1} F \not\in \overline{K'}. \hfill{\diamond}$
\medskip\\
However, we can obtain a sharper result using (3.1), (0.6) and (1.0).
\medskip\\
{\bf Theorem 3.2:}  $ B(\overline{Pic^{-1} C}) = \pi'_{*}
(\overline{Pic^{-1} C'}) \ \Longleftrightarrow \ rk(\cM/\cM^2 + t\cC)=1. $
\smallskip\\
{\bf Proof:}  $\Longrightarrow $ :  Let
$ v_0 < k_{\gamma-1} + k_1, $
and
$ \partial_{\beta} = (t^{k_{\gamma-1}} + 1)+\beta. $
Then
$ F = \cM + (t^{k_{\gamma-1}}+1)\cO, $
and
$ F $  does not contain  $ t^{-1}C, $
since
$ k_{\gamma-1} < v_0 -1 $
by the previous lemma.  If
$ \partial'_{\beta} \cO'
\longrightarrow t^{-1} F = t^{-1} \cM+I'_{\gamma-1}, $
necessarily
$ \partial_{\beta}I'_{\gamma} \longrightarrow  I'_{\gamma} $
on one hand and
$ \partial'_{\beta} I'_{j} \
\longrightarrow \ \partial'_{\beta} I'_{j-1} \forall j < \gamma -1 $
on the other hand. Therefore
$ ord_{\beta} a_0 < {\rm ord}_{\beta} a_j \forall j < v_0 - k_{\gamma-1}, $
so that
$ k_{\gamma-1} - k_{\gamma-2} \geq v_0 - k_{\gamma-1}. $
Set
$ \partial^{''}_{\beta} = t^{k_{\gamma-1} -1} + \beta, $
so that
$ \partial_{\beta} \cO \longrightarrow F^{''} = \cM + (t^{k_{\gamma-1}-1})\cO $
and
$ t^{-1}F^{''} \not \in \ Filt (\cC', \delta-1). $
If no such  $ \partial'_{\beta} $  exists,  $ F $  is already a counterexample,
so we are through.
\smallskip\\
$ \Longleftarrow $ : Let
$ v_0 = k_{\gamma-1} + k_1, \Longleftrightarrow k_i = ik_1 $
for all
$ k_i \in \Gamma $
such that
$ k_i< v_0 + 1. $
Let
$ g \in \cO $
``generate''
$  \cM, $
i.e.,
$ \cM = g\cO + t\cC. $
Let
$ \partial_{\beta} \cO \longrightarrow F, {\rm ord}_t F = h. $
If
$ h < k_1, $
necessarily
$ F \approx t^h \cO^h, $
so that  $ F $  has a deformation
$ \partial'_{\beta} \cO' \longrightarrow t^{-1} F. $
If
$ h > k_{1}-1, $
set
$ j = rk(F / t^h \cO+\cC). $
We must have
$ F = t^h \cO^h + (f_1, \cdots f_j)\cO, $
where
$ v_0 -jk_1 \leq {\rm ord}_t f_1 < {\rm ord}_t f_2
< \cdots < {\rm ord}_t f_j < v_0, \ $
and
$ v_0 - (j-i+1) k_1 \leq \ {\rm ord}_t f_i \ \forall \ i. $
We shall separate two cases: first, let
$ h = jk_1, $
so that
$ {\rm ord}_t f_i > v_0-(j-i+1)k_1 \ \forall \ i. $
Take
$ \partial'_{\beta} = t^{h-1} + \beta g^{j-\gamma}t^{-1} f_1 +
\cdots + \beta^j g^{1-\gamma}t^{-1} f_j + \beta^{j+1}, $
so that
$ \partial'_{\beta} \cO' $   gives a deformation of  $ t^{-1} F. $
If, however,
$ jk_1 < h < (j+1)k_1, $
then possibly
$ {\rm ord}_t f_i =  v_0 - (j - i +1)k_1, $
so this construction may fail. Then
$ \partial'_{\beta} \longrightarrow  t^{-1} F $
for
$ \partial'_{\beta} =
t^{-1} f_1 + \beta g^{-1} t^{-1} f_2 + \cdots
+ \beta^{j-1} g^{1-j} t^{-1} f_j + \beta^j g^{-j} t^{h-1} + \beta^{j+1}.
\hfill{\parallel} $

 From (2.5.1), (2.5.2), (2.6), (2.6.1), (2.7) and (3.2),
\medskip\\
{\bf Theorem (3.3):} The following conditions are equivalent:\\
(i) $ rk(\cM / \cM^2 + t \cC)  =  1 $ \\
(ii) $ B(\overline{Pic^{0} C})  =  \pi'_{*} \overline{Pic^{-1} C'} $ \\
(iii) $ E(\cC, \delta) = \overline{K} \ $
and
$ \ E(\cC', \delta -1) = \overline{K'}.  \hfill{\parallel}$

Next we would like to deal with the curves  $ C $  which have the
property  $ B(\overline{Pic^0 C}) = \pi'_{*}(Pic^{-1}C^{'^{=}}). $
We know that this is true if  $ C $  is locally planar, or if
$ \delta = 2. $
In general, this will be true whenever  $ C' $   is locally planar. But
\medskip\\
{\bf Lemma 3.4:}  $C'$ is Gorenstein
$\Longleftrightarrow$
either
$ v_0 > k_1 = 2 $
or
$ v_0 = k_1 < 4. $
\smallskip\\
{\bf Proof:}  Let  $ C' $   be Gorenstein, so that
$ v_1 = 2 \delta-2, v_i = v_1 - k_{i-1} $  for all  $ i > 0, $
and
$ 2 \delta -2 < v_0 < 2 \delta + 1. $
Suppose  $ v_0 > k_1. $
If
$ v_0 = 2 \delta -1, v_2 = v_0 - 1 - k_1 \in \Gamma, $
which is not possible. If
$ v_0 = 2 \delta, C $
is Gorenstein, so
$ v_1 = v_0 - k_1\Longrightarrow k_1 = 2. $
Rest of the assertion is obvious.  $ \hfill{\parallel} $
\medskip\\
These, however are not all :
\medskip\\
\underline{\bf Example 7:} \
$ \ B(\overline{Pic^0 C}) = \pi'_{*} ((Pic^{-1} C')^{=}) $
for $C$ given by \\
\centerline{\hfill
$ \cO = k \left[ [s^7, s^4 t^3, s^2 t^5, t^7] \right].
\hfill{\diamond} $}

We would like to investigate the possibility of characterizing
(or at least giving a list of) curves with this property.  This is
partially done in the following:
\medskip \\
{\bf Lemma  3.5:}  If
$  \ B(\overline{Pic^0 C}) \ = \ \pi'_{*}((Pic^{-1} C')^{=}) \  $
then either  $ C $  is locally planar or  $  \ v_0 = 2\delta -1 \  $ .
\smallskip\\
{\bf Proof:}  Let
$ F' = t^{v_0 -v_1} W', $
where
$ W' \approx \cW^0_{C', Q}, $
so that  $ F' $  is an  $ \cO' $ -module  with
$ \cC \subseteq F' \subseteq \cO^{\sim} $
and
$ rk(F'/\cC) = \delta -1. $
If  $ F \approx F' $   has a deformation by free  $ \cO $ -modules
$ \partial_{\beta} \cO, $
then
$ rk(F/\cC) = \gamma. $
On the other hand
$ rk(F/\cC) > rk(F'/ \cC)-1, $
since
$ t^{-1}\cC $  is not contained in  $ F' $  by construction. This gives
$ v_0 > 2\delta -2. $
But
$ v_0 = 2 \delta $
is the Gorenstein case, where
$ \ (Pic^0 C)^{=} \ = \ Pic^0 C \bigcup \pi'_{*} (Pic^{-1} C')^{=}, \ $
so that
$ \ (\overline{Pic^0 C}) \  $
is irreducible (by the hypothesis),
which is equivalent to  $ C $  being locally planar. $ \hfill{\parallel} $

Suppose  $ C $  is a non locally planar curve, such that
$ B(\overline{Pic^0C}) = \linebreak\pi'_{*} Pic^{-1} C'^{=} $   and
$ \cM \neq \cC $ .  We shall attempt to say as much as
possible about  $ \Gamma(\cO) $   in the following. It is not yet
possible for me to say whether or not the converse of (3.5) is
true.
$$ \begin{array}{llll}
{\rm Let} \  & F' & = & t^{v_0 - v_1} W', \ {\rm so \ that} \\
{} & \Gamma(F') & = & \{j \geq v_0 - v_1 \mid v_0 -1 -j \ {\rm
is \ not \ in} \ \Gamma \} \\
{} & {} & = & \{ v_0 - v_1, v_0 - v_2, \cdots, v_0 - k_1 : v_0 -
k_1 + 1, \cdots, v_0 -2, v_0, \cdots \}.
\end{array}
$$

\underline{\bf 3.6.} \
Suppose
$ v_0 - v_1 < k_1 $  and  $ v_1 > k_1. $
Since $F'$ has a deformation by free  $ \cO $ -modules and
$ F' \in E(\cC, \delta), v_0 - v_1 < k_1
\Longrightarrow v_0 -v_i = k_{i-1} \ \ \forall \  i > 1 $ .
This implies that if
$ \partial_{\beta} $
gives a deformation of  $ F' $
then
$ \partial_{\beta} = t^{v_0 - v_1} + \beta $
necessarily.  But then
$ v_1-1 >v_0 -k_1 -1\Longrightarrow v_1 -1 \in \Gamma(F'), $
a contradiction.
Hence either
$ v_1 = k_1 $
or
$ v_0 - v_1 = k_1. $

\underline{\bf 3.7} \
If
$ v_1 = k_1 $
then
$ \delta = k_1, \Longleftrightarrow v_0 = 2 k_1 -1. $
Set
$ F' = \cO' + t\cO' \Longrightarrow rk(F'/\cC) = \gamma+2. \  F' $
has a deformation by free  $ \cO $ -modules, so for any
$ F \approx F',  \exists \partial_{\beta} \in \cO^{\sim{*}} $
giving a deformation of  $F$ if and only if
$ F = t^2 F'. $
But then
$ F_1 \subseteq \cM \cO^{\sim} \Longrightarrow k_1 = 3. $
This is the unique case
$ \Gamma = \{ 0, 3; 5, \cdots \}. $

\underline{\bf 3.8}  If
$ v_1 > k_1, $
then
$ v_0 - v_1 = k_1. $
It can be seen easily that
$ \exists ! j* < \delta -k_1 + 1 $
such that
$ v_0 -v_i = k_i \ \ \forall \ i < j^*, \ \
v_0 - v_i = k_{i-1} \ \ \forall \ i > j^*, \ {\rm and} \ v_{j*} = \delta. $

However, all these curves do not necessarily satisfy
$ B(\overline{Pic^0 C}) = \pi'_{*} (Pic^{-1} C'^{=}) \ : $
\medskip\\
{\underline{\bf Example 8:}} If
$ \ \Gamma = \{0,k_1,k_1+2,k_1+5,\cdots 2k_1+3; \
2k_1+5,\cdots\}  \ {\rm and} \ \Gamma(F) =
\{k_1, k_1 + 1, k_1 + 3, k_1 + 5, \cdots, 2k_1 + 3; 2k_1 + 5, \cdots \} $
then
$ F \in Filt (\cC, \delta), F \not \in \overline{K}. \hfill{\diamond}$
\bigskip
\begin{center}
{\underline{\bf References}}
\end{center}
\begin{enumerate}
\item[{[AIK]}] A. Altman, A. Iarrobino and S. Kleiman,
Irreducibility of the compactified Jacobian, in ``Proceedings,
Nordic Summer School NAVF, Oslo, Aug. 5-25, 1976,'' Noordhoff,
Groningen, 1977.
\item[{[AK1]}] A. Altman and S. Kleiman, Compactifying the
Picard Scheme, Advances In Mathematics 35, 50-112 (1980).
\item[{[AK]}] A. Altman and S. Kleiman, Compactifying the Picard
Scheme, II, Amer. J. Math. 101 (1979) (``Zariski volume''),
10-41.
\item[{[D]}] C. D'Souza, ``Compactification of Generalized
Jacobians,'' Thesis, Tata Institute, Bombay, 1973.
\item[{[HK]}] Herzog and Kunz, Die Wertehalbgruppe eines Lokalen
Ringes der Dimension I, Springer Verlag, 1971.
\item[{[K]}] H. Kleppe, Picard scheme of a curve and its
compactification.  Thesis, M.I.T., 1981.
\item[{[MM]}] A. Mayor and D. Mumford, Further Comments on
Boundary Points, Amer. Math. Soc. Summer Institute, Woods Hole,
Mass., 1964.
\item[{[PS]}] G.Pfister and J.H.M. Steenbrink, Reduced Hilbert Scheme for
Irreducible Curve Singularities, to appear, Journal of Pure and Applied
Algebra.
\item[{[R]}] C. Rego, The Compactified Jacobian, Ann. Scient.
Ec. Norm. Sup., 4 $ ^4 $  serie, t. 13, 1980, p. 211 a 224.
\item[{[SK]}] S. Kleiman: The structure of the compactified
Jacobian : a review and an announcement; seminari di geometria,
1982-83; Universita delgi Studi di Bologna, Dipartmento
Mathematica, 1984.
\end{enumerate}
\end{document}